\documentclass[11pt]{article}
\usepackage{amsmath,amssymb,amsfonts,latexsym,lscape,rawfonts}
\numberwithin{equation}{section} 
\usepackage{lscape}

\newcommand{\beqn}{\begin{eqnarray}}
\newcommand{\eeqn}{\end{eqnarray}}

\newcommand{\beq}{\begin{equation}}
\newcommand{\eeq}{\end{equation}}

\newcommand{\R}{\mathbb{R}}

\newcommand{\Cut}{\mathrm{Cut}}
\newcommand{\Inj}{\mathrm{Inj}}
\newcommand{\HH}{\mathbb{H}}
\newcommand{\C}{\mathbb{C}}

\renewcommand{\baselinestretch}{1.0}
\newtheorem{thm}{Theorem}[section]
\newtheorem{prop}[thm]{Proposition}
\newtheorem{lemma}{Lemma}[section]

\newtheorem{cor}[thm]{Corollary}

\newtheorem{definition}[thm]{Definition}
\newcommand{\E}{\text{Exp}}

\begin{document}
\title {{\bf{\Large An intrinsic proof of Gromoll-Grove diameter rigidity theorem}}
}

\vspace{.5cm}
\author{Jianguo Cao\footnote{The first author is supported in part by an NSF
grant. Mailing Address: Department of Mathematics, University of
Notre Dame, Notre Dame, IN 46556, USA; } \ and Hongyan
Tang\footnote{Mailing Address: Department of Mathematical Sciences,
Tsinghua University, Beijing 100084, P. R. China. }
 }

\date{}
\maketitle

\begin{center}
{\it Dedicated to Professor Karsten Grove on his sixtieth birthday}
\end{center}


\vspace{.3cm}

{\small ~~~{\bf Key words and phrases}:   positive curvature,
diameter, simply-connected, rigidity

~~~{\bf MSC 2000}: } 53C20
\medskip\medskip

\section{Introduction}
We will present a new proof of the following Gromoll-Grove
diameter rigidity theorem.

\medskip
\noindent{\bf Theorem A} \ {\it Let $M^n$ be a simply connected
Riemannian manifold with sectional curvature $K \geq 1$. Suppose
that $\mathrm{Diam}(M^n)=\frac{\pi}{2}$ and $M^n$ is not
homeomorphic to a sphere $S^n$. Then $M^n$ is isometric to one of
$\C P^{\frac{n}{2}}$, $\HH P^{\frac{n}{4}}$ or $\C aP^2$, i.e.,
$M^n$ is isometric to a projective symmetric space over complex
numbers, or quaternion numbers or Calay numbers. }

\medskip

Our proof  does not use any loop spaces, which is totally
different from [Wil]. Among other things, we use  the Hessian comparison theorem for distance
functions and the spherical metric on the tangent space instead, see Section 3 below.
Although our new proof is longer than its earlier version,
 the most of arguments below remain to be elementary and
 self-contained.

\section{The Gromoll-Grove fibration}

\quad We need to recall some known results from [GG1], in order to complete
the proof of Theorem A. The results of [GG1] are related to the following example.

\bigskip
\noindent {\bf Example 2.0.} (1) Let $M^n = \C P^{\frac{n}{2}}$ with
the classical Fubini-Study metric and diameter $\frac{\pi}{2}$. Let
$B_r(p)$ be the metric ball of radius $r$ and center $p$ in $\C
P^{\frac{n}{2}}$, and let $S_r(p) = \partial B_r(p)$ be the metric
sphere of radius $r$ centered at $p$. It is well-known that
$S_{\frac{\pi}{2}}(p)$ is isometric to $\C P^{\frac{n}{2} -1}$.

For each $p \in  \C P^{\frac{n}{2}}$, we consider
a polar coordinate $\{ (r, \Theta) \}$ of the tangent space $T_p(\C P^{\frac{n}{2}})$ and the exponential map
$\E_p: T_p(\C P^{\frac{n}{2}}) \to \C P^{\frac{n}{2}}$.

    Let us choose a spherical metric
    $$
   g_1 = dr^2 + (\sin r)^2 d \Theta^2
    $$
on the set $B_{\pi}(0) = \{ ( r, \Theta) | 0 \le r < \pi, \Theta \in
S^{n-1}\}$, where $d \Theta^2$ is the canonical metric of constant
curvature $1$ on the unit sphere $S^{n-1}$. With respect to the
spherical metric $g_1$, the exponential map
$$
\begin{array}{rrcl}
\E_p : & S_{\frac{\pi}{2}}(0) & \rightarrow & S_{\frac{\pi}{2}}(p) \\
        &\frac{\pi}{2} \Theta   & \rightarrow &
\mathrm{Exp}_p(\frac{\pi}{2} \Theta)
\end{array}
$$
is a Hopf fibration.

Furthermore,  for each $q \in  S_{\frac{\pi}{2}}(p)$, the fiber
$\E_p^{-1}(q)$ is a great circle in the equator
$(S_{\frac{\pi}{2}}(0), g_1)$ of the unit sphere $S^n =
(\bar{B}_\pi(0), g_1)$.

(2) We are going to elaborate the above construction by replacing
the point $p$ by a totally geodesic submanifold $\C P^m \subset \C
P^{\frac{n}{2}}$ with $1 \le m < \frac{n}{2} - 1$, for the case
$\frac{n}{2} \ge 3$. We let $U_r(\C P^m) = \{ z \in \C
P^{\frac{n}{2}} \, | \, d(z, \C P^m) < r \}$ be the tubular
neighborhood and $\partial [U_r(\C P^m)]$ its boundary.

Then $\partial [U_{\frac{\pi}{2} }(\C P^m)]$ is isometric to a
totally geodesic $\C P^{m'} \subset \C P^{\frac{n}{2}}$ with $m' =
\frac{n}{2} - m - 1$. In this case, for each pair $p \in \C P^m$ and
$q \in \C P^{m'}$ with distance $d(p, q) = \frac{\pi}{2}$, we still
have that the fiber $\E_p^{-1}(q)$ is a great circle in the equator
$(S_{\frac{\pi}{2}}(0), g_1)$ of the unit sphere $S^n = (B_\pi(0),
g_1)$, where $\bar{B}_r(0) \subset T_p(\C P^{\frac{n}{2}}  )$.

In fact, $\C^{\frac{n}{2} +1 } $ has a decomposition
$\C^{\frac{n}{2} +1 } = \C^{m+1} \times \C^{m' + 1}$. Such a
decomposition induces a spherical join of $S^{2m+1}$ and
$S^{2m'+1}$. More precisely, for each unit vector $\vec u \in
S^{n+1} \subset \C^{\frac{n}{2} +1 }$, there are $\vec v \in
S^{2m+1}$ and    $\vec w \in S^{2m'+1}$
$$
\vec u = (\cos r) \vec v + (\sin r) \vec w
$$
for some $r \in [0, \frac{\pi}{2}]$. One can write $S^{n+1} =
S^{2m+1} \star S^{2m'+1}$, where $ \frac n2 = m + m'+1$. It follows
that $\C P^{\frac{n}{2}}$ can be viewed as the {\it ``projective
join" } of $\C P^m$ and $\C P^{m'}$.

The pair of sub-manifolds $\{\C P^m, \C P^{m'} \}$ with $d( \C P^m,
\C P^{m'}) = \frac{\pi}{2}$ above is called a {\it dual pair} of
convex subsets of $\C P^{\frac n2}$ in [GG1].

(3) When $M^n$ is isometric to either $\HH P^{\frac{n}{4}}$ or $\C
aP^2$, there are similar decompositions. \hfill{Q.E.D.}

\medskip
\bigskip

Inspired by Example 2.0, we consider the convexity of subset $[M -
B_r(p)]$, without  the  assumption $\mathrm{Diam}(M) =
\frac{\pi}{2}$. Let $\mathrm{Inj}_M(x)$ denote the injectivity
radius of $M$ at $x$.

\begin{prop}  Let $M$ be a complete smooth Riemannian manifold with sectional
curvature $\ge 1$ and $\mathrm{Diam}(M) \ge\frac{\pi}{2}$. Suppose
that $\sigma: [0, \ell] \to M$ is a length-minimizing geodesic of
unit speed from $x$. Then, for any $0< r < \ell$, the second
fundamental form of $S_r(x)$ at $\sigma(r)$ with respect to the
normal vector $\sigma'(r)$ is less than or equal to $ \cot (r) I$ at
$\sigma(r)$ in the barrier sense, where $I$ is the identity matrix.

Consequently, if $Inj_M(x) \ge \frac{\pi}{2}$,  then $[M -
B_{\frac{\pi}{2}}(x) ]$ is a convex subset of $M$. In addition, if
$\mathrm{Diam}(M) \ge \ell > \frac{\pi}{2} $, then $[M - B_\ell(x)
]$ is strictly convex.
\end{prop}
{\bf Proof.} This is a direct consequence of the Hessian comparison
(see [Pe, p145]) for the distance function.      \hfill Q.E.D.

\begin{prop}\label{thm:00} ([GG1])
Let $M$ be a complete smooth Riemannian manifold with sectional
curvature $\ge 1$ and $\mathrm{Diam}(M) = \frac{\pi}{2}$. If $M$ is
simply-connected and if $M$ has the integral cohomology ring of
either $\C P^2$, $\HH P^2$ or the Cayley plane $\C aP^2$, then there
exists at least one point $p \in M$ with injectivity radius
$\Inj_M(p) \ge \frac{\pi}{2}$.
\end{prop}

When $d(p, q)=\mathrm{Diam}(M) = \frac{\pi}{2}$, the subset
$S_{\frac{\pi}{2}}(p)$ is a critical sub-manifold of the distance
function $f(x) = d(x, p)$. If $\Inj_M(p) \ge \frac{\pi}{2}$,  by
Proposition \ref{thm:00} above, $S_{\frac{\pi}{2}}(p)$ is a totally
geodesic submanifold. In fact, the dual convex subset
$S_{\frac{\pi}{2}}(p)$ has some extra properties (cf. [GG1]), which
we recall in the sequel.

Following [GG1], for $A \subset M$ we let
$$
A' = \{y\in M \ |\ d(y, A) = \frac{\pi}{2}\}.
$$
The following result was  also stated in [GG1].

\begin{prop}\label{thm:01} ([GG1])
Let $M$  be a simply connected Riemannian manifold with sectional
curvature $\geq 1$ and $\mathrm{Diam}(M)=\frac{\pi}{2}$. Suppose
that the injectivity radius $\Inj_M(p)$ of $M$ at $p$ is equal to
$\frac{\pi}{2}$ and that $M$ is not homeomorphic to a sphere. Then

(1) $M$ has integral cohomology ring of either $\C P^{\frac n2}$,
$\HH P^{\frac n4}$ or the Cayley plane $\C aP^2$;

(2) if $A=\{p\}$, then $A'=\{y \ | \ d(y, p) = \frac{\pi}{2}\}$ is a
closed totally geodesic submanifold of positive dimension;

(3) if $A=\{p\}$, then $(A')'= A$ and the cut radius $\Cut_{M}(A')$
is equal to $\frac{\pi}{2}$ as well;

(4)  if $S_p (M)= \{\vec v \in T_p M\ | \ \|\vec v\|=1\}$ then
$$
\begin{array}{rrcl}
\pi_p = \widetilde{\mathrm{Exp}}_p: & S_p(M) & \rightarrow & A' \\
                          & \vec v   & \rightarrow & \mathrm{Exp}_p(\frac{\pi}{2}\vec v)
\end{array}
$$
is a Riemannian submersion;

(5) if the Riemannian submersion $\pi_p : S_p(M) \rightarrow  A' $
is a great circle fibration, then $M$ is isometric to either $\C
P^{\frac n2}$, $\HH P^{\frac n4}$ or  the Cayley plane $\C aP^2$.
\end{prop}

\bigskip

Notice that $ \{ q \}''' = \{ q \}'$ holds for all $q \in M$, when
$\mathrm{Diam}(M)=\frac{\pi}{2}$. We choose $A' =  \{ q \}'$ and $A
= A''$. It is possible that $\min \{ \dim A, \dim A' \} > 0$, see
Example 2.0 above. If $M$ is allowed to be non-simply-connected, and
if $M^3$ is a lens space, then $\min \{ \dim A, \dim A'\} = 1$. In
other words, it might be difficult to find a point $p$ with
$\Inj_M(p) = \frac{\pi}{2}$ when $\mathrm{Diam}(M) = \frac{\pi}{2}$.
We can not choose $A$ with $\dim A = 0$ at the first place.

Thus, we need to describe the remaining case of $\min\{\dim A, \dim
A'\} > 0$, where $A'' = A$ and $\{A, A'\}$ is a pair of dual convex
subsets. It was shown in [GG1] that both $A$ and $A'$ are connected
totally geodesic submanifolds without boundaries.

In what follows, we always let
$$
S^\perp_q(B, M) = \{ \vec v \in T_p M \, \, | \,  \vec v  \perp
T_q(B), | \vec v | = 1 \}
$$
be the unit normal bundle of $B$ in $M$, when $A'$ is a submanifold
of M.

\begin{prop}\label{thm:02} ([GG1])
Let $M$ be a simply-connected Riemannian manifold with sectional
curvature $\ge 1$ and diameter $\mathrm{Diam}(M^n)=\frac{\pi}{2}$.
For any $z \in M$ with $S_{\frac{\pi}{2}}(z) \neq \emptyset$, we let
$A' = S_{\frac{\pi}{2}}(z)$ and $A = A''$. Suppose that $M^n$ is not
homeomorphic to $S^n$. Then

(1) both $A$ and $A'$ are simply-connected;

(2) the cut radius $\mathrm{Cut}_{M}(A)$ of $A$ in $M$ is equal to
$\frac{\pi}{2}$ and $\mathrm{Cut}(A) = A'$;

(3) the cut radius $\mathrm{Cut}_{M}(A')$ of $A'$ in $M$ is equal to
$\frac{\pi}{2}$ and $\mathrm{Cut}(A') = A$;

(4) if  $\dim(A') > 0$, then
$$
\begin{array}{rrcl}
\pi_p = \widetilde{\mathrm{Exp}}_p: & S_p^\perp(A, M) & \rightarrow & A' \\
                          & \vec v   & \rightarrow & \mathrm{Exp}_p(\frac{\pi}{2}\vec v)
\end{array}
$$
is a Riemannian submersion; similarly, if $\dim A > 0$ then $\pi_q :
S_q^\perp(A', M) \rightarrow  A$ is a  Riemannian submersion for all
$q \in A'$; furthermore, $\dim[\pi_p^{-1}(q)]$ is equal to one of
$\{ 1, 3, 7\}$; $\dim M$, $\dim A$ and $\dim A'$ are even integers;

(5) if the Riemannian submersion $\pi_p : S_p^\perp(A, M)
\rightarrow A' $ with $\dim A' > 0$ is a great circle fibration for
all $p \in A$ and if the Riemannian submersion $\pi_q :
S_q^\perp(A', M) \rightarrow  A$ is a great circle fibration for all
$q \in A'$ whenever $\dim A > 0$, then $M^n$ is isometric to one of
symmetric spaces $\C P^{\frac n2}$, $\HH P^{\frac n4}$ or $\C aP^2$.
\end{prop}

\begin{definition}
 Let $M^n$ be a simply connected
Riemannian manifold with sectional curvature $\geq 1$. Suppose that
$\mathrm{Diam}(M^n)=\frac{\pi}{2}$,  $A'' = A$ and  $(p, q) \in A
\times A'$. When $\dim A' > 0$,  the Riemannian submersion
$$
\begin{array}{rrcl}
\pi_p : & S_p^\perp(A, M) & \rightarrow & A' \\
        & \vec v   & \rightarrow &
\mathrm{Exp}_p(\frac{\pi}{2}\vec v)
\end{array}
$$
is called the Gromoll-Grove fibration with the total space
$S_p^\perp(A, M)$.

 Similarly, when $\dim A > 0$, the fibration $\pi_q : S_q^\perp(A', M)
\rightarrow  A $ is called the Gromoll-Grove fibration as well.
\end{definition}

In next section, we will show that the Gromoll-Grove fibration
$$
S^k \rightarrow S_p^\perp(A, M) \rightarrow A'
$$
is a great circle fibration for some $k \in \{ 1, 3, 7\}$ whenever
$\dim(A') > 0$; and hence $M$ must be isometric to a symmetric space
by [Ran].

\section{The Gromoll-Grove fibration is isometrically congruent to
a Hopf fibration}

In this section, we will use a new method to show that the
Gromoll-Grove fibration is isometrically congruent to a great circle
fibration.

Throughout this section, the origin of $T_p M \approx \R^n$ is
denoted by $0_p$. We will always use a spherical metric $g_1$ on a
ball $B_\pi(0_p) \subset T_p M$:
$$
g_1 = dr^2 + (\sin r)^2 d \Theta^2
$$
where $\{(r, \Theta) \}$ is the polar coordinate system of $T_pM
\approx \R^n$.

We consider the possibly tear-drop shaped fibres in the manifold
$M$, see Section 2 above. For each pair of $p \in A$ and $q \in A'$,
we let
$$
\Sigma_{p, q} = \{ \text{Exp}_p(t \vec v) \, |  \, \vec v \in
\pi_p^{-1}(q), 0 \le t \le \frac{\pi}{2}  \}.
$$
and
$$
\tilde \Sigma_{p, q} = \{ \vec w \in \E_p^{-1}(\Sigma_{p, q}) \, |
\, \| \vec w\|  \le \frac{\pi}{2}\}
$$
 be the truncated tangential cone of $\Sigma_{p, q}$ at $p$.

 Our goal is to show that $\tilde \Sigma_{p, q}$ is totally geodesic in
 $(B^\perp_\pi(0_p), g_1) \subset S^n $ and hence $\partial [\tilde \Sigma_{p, q} ]$
 is totally geodesic in $S^n$. Consequently, $\pi^{-1}_p(q) $ is a $k$-dimensional circle in
$S^{n-1}$, where $k$ is one of $\{1, 3, 7\}$.

\medskip

There are three elementary steps to show that $\pi^{-1}_p(q) $ is a
$k$-dimensional circle in $S^{n-1}$.

\medskip
\noindent
 {\it Step 1.} We will show that {\it ``if $\tilde \Sigma_{p, q}$ has the first focal radius
 $\ge \frac{\pi}{2}$ in $ S^n$ at all $ z \in \tilde \Sigma_{p, q} $ with $0< |z| < \frac{\pi}{2}$,
 then $\tilde \Sigma_{p, q}$ is a smooth totally geodesic submanifold of $S^n $".}

\medskip
\noindent {\it Step 2.}  We will make the following elementary
observation. {\it Suppose contrary, $\tilde \Sigma_{p, q}$ had the
first focal radius $ 0 < t_0< \frac{\pi}{2}$ in $ (B^\perp_\pi(0_p),
g_1) \subset S^n$ at some $ z \in  \tilde \Sigma_{p, q} $ with $0 <
|z| < \frac{\pi}{2} $. Then there would be a Jacobi field $\{ J(t)
\}$ along a normal geodesic $\sigma_{z,\vec h }(t) = \E^{S^n}_z(t
\vec h)$ such that $\vec h \perp T_z(\tilde \Sigma_{p, q}  )$,
$|\vec h| = 1$ and $J'(0)  \in T_z(\tilde \Sigma_{p, q}  )$.}

Thus, we consider a special class of Jacobi fields with  {\it extra}
initial conditions on $J'(0)$:
      $$
\Gamma_{\sigma_{z, \vec h}, \tilde \Sigma_{p, q}} = \{ J \, | J'' +
R(\sigma', J) \sigma' = 0, \langle J'(0), X \rangle = - \langle \vec
h, \nabla_{J(0) } X \rangle, \mathrm{ for } \, \, \mathrm{ all } \,
X \in T_z(\tilde \Sigma_{p, q})\}
      $$
 and
\beqn\label{eq:20} \Gamma^0_{\sigma_{z, \vec h}, \tilde \Sigma_{p,
q}} = \{ J \in  \Gamma_{\sigma_{z, \vec h}, \tilde \Sigma_{p, q}}, |
 J'(0) \in T_z(\tilde \Sigma_{p, q})
\}.
\eeqn
 It will be shown
 $$\dim[\Gamma^0_{\sigma_{z, \vec h}, \tilde \Sigma_{p, q}}  ] =
      \dim [ \tilde \Sigma ] = k+1.
      $$
       This step is applicable to all $(k+1)$-dimensional submanifold
      $\tilde \Sigma \subset S^n$, which is elementary.

\medskip
\noindent
 {\it Step 3.} In this final step, we use Hessian comparison theorem to show that,
 {\it ``if $\pi_p: S^\perp_p(A, M) \to A'$ is a Riemannian submersion, then,  for all non-trivial Jacobi field
 $J \in \Gamma^0_{\sigma_{z, \vec h}, \tilde \Sigma_{p, q}}$, we have $J(t) \neq 0$ for all
 $t \in (0, \frac{\pi}{2})$.''} It follows that  $\tilde \Sigma_{p, q}$ has the first focal radius
 $\ge \frac{\pi}{2}$ and hence totally geodesic in $S^n$. This completes the proof of Grove-Gromoll
 diameter rigidity Theorem.

\bigskip
\bigskip

Here are the details for each step.

\bigskip

\noindent {\bf Step 1.} We present a sufficient condition for
totally geodesic property.

A subset $ C \subset M$ is $a$-convex in $M$ if, for all geodesic
segments $\sigma: [0, \ell] \to M$ of length $\ell < a$ with
endpoints in $C$, one has $\sigma([0, \ell]) \subset C$.

\begin{prop}\label{thm:03} If $\tilde \Sigma_{p, q}$ has the first focal radius $\ge \frac{\pi}{2}$ in
$S^n = (B_\pi(0_p), g_1) $ at all $ z \in  \tilde \Sigma_{p, q} $
with $0 < |z| < \frac{\pi}{2} $ where $B_\pi(0_p) \subset T_p(M^n)$,
then

 (1) $\tilde \Sigma_{p, q}$ is a smooth totally geodesic submanifold
with boundary in $S^n = (B_\pi(0_p), g_1)$; Moreover, $\pi_p^{-1}(q)
\approx [\partial  \tilde \Sigma_{p, q}]$ is a totally geodesic
great $k$-dimensional circle in $S^n$.

 (2) The injectivity radius $\Inj_{M}(q)$ of  $q$ in $M$ is equal to $\frac{\pi}{2}  $
 for all $q \in A'$.
\end{prop}
{\bf Proof.} (1) Let $\vec h_0 \in S^\perp(\tilde  \Sigma_{p, q},
S^n)$ be a unit norm vector of $ \tilde \Sigma_{p, q}$ at $z_0$
and $\sigma_0(t) = \E_{ z_0}( t \vec h_0 )$. Let $\E: S^\perp(
\tilde \Sigma_{p, q}, S^n ) \times [0, \infty) \to S^n$ be the
exponential map along the normal bundle near $(z_0, \vec h_0; t)$
for $ t \ge 0$. Suppose that $ \zeta: ( -\delta, \delta)
\to S^\perp( \tilde  \Sigma_{p, q}, S^n ) $ is a curve with
$\zeta(0)= ( z_0, \vec h_0)$ and $\zeta(s)= (z(s), \vec h(s))$.
Then $F(t, s) = \E_{z(s)}[t \vec h(s)]$ gives rise to a Jacobi
field $\{ J(t) \}$ defined by
$$
J(t) = \frac{\partial F}{\partial s}(t, 0)
$$
along $\sigma_0$.

Our goal is to show that, under our assumption, we have

\beqn\label{eq:5}
 \langle
J(0), \nabla_{J(0)} \vec h(s) \rangle = \langle J(0), J'(0)
\rangle \ge 0. \eeqn

Since $\tilde \Sigma_{p, q}   $ is not a hypersurface, we consider a tubular neighborhood of
 $\tilde \Sigma_{p, q}   $. Choose $\varepsilon_1 $ sufficiently small so that
$B_{\varepsilon_1  }( z_0) \cap \tilde \Sigma_{p, q}$ is an embedded
$(k+1)$-dimensional ball. Let $\varepsilon_0$ be the cut-radius of
$ \tilde \Sigma_{p, q}\cap B_{\varepsilon_1  }(z_0)  $. Choose $ \varepsilon < \frac 18 \min\{\varepsilon_0,
\varepsilon_1 \}$. Then there is a nearest point projection from $B_\varepsilon(z_0) \to \tilde \Sigma_{p, q}$.
If $\{ G(., s)\}$ is an 1-family variation of $\sigma_0|_{[\varepsilon, \ell]}$, which are orthogonal to
$\partial [U_\varepsilon( \tilde \Sigma_{p, q} )    ]$, then by using the nearest point projection, such a family
$\{G(., s)\}$ can be extended as an 1-family of normal geodesics $\{ F(., s) \}$ from $\tilde \Sigma_{p, q}$
with $F(., 0) = \sigma_0(.)$.

Hence we see that the hypersurface $\partial [U_\varepsilon( \tilde
\Sigma_{p, q} )    ]$ has focal radius  $\ge \frac{\pi}{2}  -
\varepsilon$ along $\sigma_0$, where $U_\varepsilon( C) = \{ y \in
S^n \, | \,  d(y, C) < \varepsilon\}$.

To prove \eqref{eq:5}, it is sufficient to
\beqn\label{eq:6}
 \frac{ \langle J( \varepsilon ), J'( \varepsilon)
\rangle}{|J( \varepsilon )   |^2} \ge  - \tan  \varepsilon. \eeqn

Let $\lambda(\varepsilon) $ be an eigenvalues of $II_\varepsilon(X,
Y) = - \langle \nabla_XY , - \sigma_0'(\varepsilon)\rangle$. For
\eqref{eq:6}, it is sufficient to show that $\lambda(\varepsilon)
\ge  - \tan  \varepsilon$.

We may isometrically embed $S^n$ into $\R^{n+1}$ as
$S^n \approx \{ \vec w \, | \vec w \in \R^{n+1}, \, | \vec w | = 1 \}$ and identify $0_p \in T_p(M^n)$ with
the North pole $e_{n+1} = (0, ..., 0, 1) $.
In the unit sphere $S^n \subset \R^{n+1}$, any geodesic of unit speed can
be written as $\sigma(t) = (\cos t) \sigma(0) + (\sin t) \sigma'(0)$.
Thus, any Jacobi field can be expressed as
$ J(t) = (\cos t) J(0) + (\sin t) J'(0)$. If $J'( \varepsilon) = \lambda(\varepsilon) J( \varepsilon )$,
then we have $J(t) = [\cos (t- \varepsilon) + \lambda(\varepsilon) \sin ( t- \varepsilon )  ] J( \varepsilon)$.
The equality $J(t_0) = 0$ holds if and only if $t_0 = \cot^{-1}[- \lambda(\varepsilon)] + \varepsilon$.
By our assumption, $t_0 \ge \frac{\pi}{2}$. It follows that
$$
- \lambda(\varepsilon) = \cot [t_0 - \varepsilon  ]\le \cot
[\frac{\pi}{2} -  \varepsilon] = \tan  \varepsilon
$$
This completes the proof of \eqref{eq:5} and \eqref{eq:6}.

Hence, we showed that $\tilde \Sigma_{p, q}$ is totally geodesic at
$z$ with $0< d(z, 0_p) < \frac{\pi}{2}$. It remains to show that
$\partial [\tilde \Sigma_{p, q}]$ is a $k$-dimensional great circle
in $S^n$.

For this purpose, we let $S^{n-1} = \{ (\vec v, 0) \,  | \,  (\vec
v, 0) \in S^n \subset \R^{n+1} \}$ be the equator of $S^n$. Let
$\Psi: S^n - \{\pm e_{n+1}\} \to S^{n-1}$ be the nearest point
projection to the equator $S^{n-1} $ given by $ \Psi(z) = \frac{ z -
\langle z, e_{n+1} \rangle e_{n+1} } {|z - \langle z, e_{n+1}
\rangle e_{n+1}|  }$. It is easy to see that $\Psi $ takes a
geodesic segment in $S^n$ to a arc of a great circle in $S^{n-1}$.
Since $\tilde \Sigma_{p, q}$ is totally geodesic at $z$ with $0<
d(z, 0_p) < \frac{\pi}{2}$, $\Psi(\tilde \Sigma_{p, q}  )$ must be
contained in a totaly geodesic subset in the equator $S^{n-1}$.
However, it is easy to see that $ \Psi(\tilde \Sigma_{p, q}  )
=\partial [ \tilde \Sigma_{p, q}  ]$. It follows that $\partial [
\tilde \Sigma_{p, q}  ]$ is totally geodesic in $S^n$. Consequently,
$\pi_p^{-1}(q)$ is a great circle in $S^{n-1} \subset \R^n$.

\medskip

(2) We first observe $\Inj_{A'}(q) = \frac{\pi}{2}$, due to [Ran]. Here is a direct proof of
$\Inj_{A'}(q) = \frac{\pi}{2}$ without using results of [Ran].

We now consider the Riemannian submersion $\pi_p: S^\perp_p(A, M   )
\to A'$, where $S^{n-1} =(\partial B_{\frac{\pi}{2}}(0), g_1)$ is
the equator of $S^n \subset \R^{n+1}$. Since $\pi^{-1}(q)$ is
totally geodesic, it is a great circle. We still isometrically embed
$S^n$ into $\R^{n+1}$ as above. It follows that the linear sub-space
$\mathrm{Span}\{\pi^{-1}(q)\}$ spanned by $\pi^{-1}(q)$ is isometric
to a $(k+1)$-dimensional $\R^{k+1}_q$.

For each geodesic segment of unit speed $\hat \sigma: [0, \frac{\pi}{2}] \to A'$ from $q$ to $y = \hat
\sigma( \frac{\pi}{2} )$, we will show that $d_{A'}(q, y) = \frac{\pi}{2}$.

Let $\tilde \sigma: [0, \frac{\pi}{2}] \to S^{n-1} $ be a horizontal
lift of $\hat \sigma$. As we pointed out above, we can write $\tilde
\sigma(t) = \cos t \tilde q + \sin t \tilde \sigma'(0)$, where
$\tilde \sigma'(0) = \tilde y \perp \tilde q$. At time $t
=\frac{\pi}{2}$, the vector $\tilde \sigma'(\frac{\pi}{2})$ becomes
horizontal. Thus, $\tilde q = \tilde \sigma'(\frac{\pi}{2})$ is
orthogonal to $T_{\tilde y}(\pi_p^{-1}(y)) \subset R_y^{k+1}$, where
$\tilde y = \tilde \sigma(\frac{\pi}{2}) \in \pi_p^{-1}(y)$. Recall
that $\tilde y = \tilde \sigma'(0) \perp \tilde q$. Hence,  $\tilde
q \perp \R^{k+1}_y$.

Suppose contrary, if $d_{A'}(q, y) = \alpha < \frac{\pi}{2}$. Then
there would be another length-minimizing geodesic $ \hat \sigma_2:
[0, \alpha] \to A'$ from $q$ to $y$. Using the horizontal lift
$\tilde \sigma_2$ of $\hat \sigma_2$ with the initial point $\tilde
q$, we would be able to find $\tilde z = \tilde \sigma_2( \alpha)
\in \pi_p^{-1}(y)$. It would follow that the angle between $\tilde
q$ and $\tilde z$ is equal to $ \alpha <  \frac{\pi}{2}$, which
contradicts to the fact $\tilde q \perp \R^{k+1}_y$. Thus, any
geodesic segment $\hat \sigma: [0, \frac{\pi}{2}] \to A'$ of unit
speed is length-minimizing, and hence $\Inj_{A'}(q) =
\frac{\pi}{2}$.

\medskip

Let us now further prove $\Inj_M(q) = \frac{\pi}{2}$. Let  $\sigma:
[0, \frac{\pi}{2}] \to M$ be any geodesic segment with $\sigma(0) =
q$ and unit speed. If $\sigma'(0) \perp A'$, then by Proposition
\ref{thm:02}, $z=\sigma( \frac{\pi}{2} ) \in A$ and hence $d(q,
\sigma( \frac{\pi}{2} )) =  \frac{\pi}{2} $. Thus, $\sigma: [0,
\frac{\pi}{2}] \to M$ is length-minimizing in this case.

If $\sigma'(0) = (\cos \beta) \vec v + (\sin \beta) \vec h$ for some
$\vec v \perp A'$,  $\vec h \in T_q(A')$ and $0 < \beta<
\frac{\pi}{2}  $, we let $\Psi_A: [M-A'] \to A$ be the nearest point
projection, and let $\Psi_{A'}: [M-A] \to A'$ the nearest point
projection. Let
$$
z = \sigma( \frac{\pi}{2} ).
$$
For $y = \Psi_{A'}(z) = \Psi_{A'}(\sigma( \frac{\pi}{2} )  )$, by
Lemma 3.1 of [GG1], $\{q, y, z\}$ and $\Psi_{A'}(\sigma( \R))$ are
contained in a totally geodesic 2-sphere. Moreover,
$\Psi_{A'}(\sigma( [0, \frac{\pi}{2}] ))$ is a geodesic segment of
length $\frac{\pi}{2}$. Thus, since the injectivity radius of $A'$
is equal to $\frac{\pi}{2}$, one has that $d_{A'}(y, q)$ is equal to
the length of $\Psi_{A'}(\sigma( [0, \frac{\pi}{2}] ))$, which is
$\frac{\pi}{2} $. Let $x = \Psi_A(z) \in A$. It is clear that $d(x,
q) \ge d(A, q) = \frac{\pi}{2}$. Hence, we have $\{x, y\} \subset
[\partial B_{\frac{\pi}{2}}(q)] = \{ q \}'$.

We already showed that  $\{x, y\}\subset \{q\}'$ holds. It now
follows from Proposition 1.3 of [GG1] that $ \{ q \}' $ is
$\pi$-convex.  Because $z$ lies on a geodesic segment of length
$\frac{\pi}{2} < \pi$ from $x$ to $y$ and $\{x, y\} \subset \{q\}'$,
by the $\pi$-convexity of $\{q\}'$ we obtain that $z \in \{ q \}'$.
Therefore, any geodesic segment $\sigma: [0, \frac{\pi}{2}] \to M$
of unit speed from $q$ is length-minimizing for all cases. The
assertion of $\Inj_M(q) = \frac{\pi}{2}$ is proved. \hfill{Q.E.D.}

\bigskip
\noindent {\bf Step 2.} For the convenience to the reader, we
include the detailed proof of the following elementary result.

\begin{prop} Let $\tilde \Sigma \subset S^n$ be a $(k+1)$-dimensional submanifold which is smooth at
 $z\in \tilde \Sigma$. Suppose that $\vec h \in T_z(S^n)$ is a unit normal vector of $\tilde \Sigma$ at $z$,
 $\sigma_{z, \vec h} (t) = \E_z^{S^n}(t \vec h)$ and let $\Gamma^0_{\sigma_{z, \vec h}, \tilde \Sigma_{p, q}}$
 be as in \eqref{eq:20} above. Then

(1) $\dim [\Gamma^0_{\sigma_{z, \vec h}, \tilde \Sigma_{p, q}}] = k +1 $;

(2) If $\tilde \Sigma$ has the first focal radius $t_0 <
\frac{\pi}{2}$ along $\sigma_{z, \vec h}$, then there must be a
non-trivial Jacobi field $\{J(t) \}$ along $\sigma_{z, \vec h}$ with
$J'(0) = (- \cot t_0) J(0) \in T_z( \tilde \Sigma)$, and hence $J
\in \Gamma^0_{\sigma_{z, \vec h}, \tilde \Sigma_{p, q}}$.
\end{prop}
{\bf Proof.} (1) Let $N(\tilde \Sigma)|_{B_\varepsilon(z)} = \{ (y,
\vec w)  | y \in \tilde \Sigma, d(y, z) < \varepsilon, \vec w \perp
T_y( \tilde \Sigma )\}$ be the normal bundle of $\tilde \Sigma$ near
$z$. Suppose that $G=\E^{S^n}: N(\tilde \Sigma)|_{B_\varepsilon(z)}
\to S^n$ is the exponential map of $S^n$ restricted to the normal
bundle of $\tilde \Sigma$. For any curve $\zeta: (-\delta, \delta  )
\to N(\tilde \Sigma)$ with $\zeta(0) = (z, \vec h)$ with $\zeta(s) =
(y(s), \vec w(s))$, there is an 1-family of geodesics given by $F(t,
s) = G( y(s), t\vec w(s)  ) = \E_{y(s)}[t \vec w(s)]$.

Let  $J(t) = \frac{\partial F }{\partial s}(t, 0) = G_*|_{(z, t \vec
h)} \zeta'(0)$. Since $\frac{\partial F }{\partial t}(0, s) = \vec
w(s) \perp \tilde \Sigma$ for all $s \in (-\delta, \delta  )$, by
the Gauss-Coddazi equation we obtain that, for all $X \in T_z(
\tilde \Sigma)$,
$$
\langle J'(0), X \rangle = \langle \vec w'(0), X \rangle = - \langle \vec w(0), \nabla_{y'(0)}X \rangle
= - \langle \vec h, \nabla_{J(0)}X \rangle
$$
holds. Hence, the tangential component of $J'(0)$ is uniquely determined by the second fundamental
form of $\tilde \Sigma$:
\beqn\label{eq:11}
\langle J'(0), X \rangle = - II_{\vec h}(J(0), X) = - \langle \vec h, \nabla_{J(0)}X \rangle
\eeqn
for all $X \in T_z( \tilde \Sigma)$. Let us consider the classical Weingarten map $W^{\vec h}:
T_z(\tilde \Sigma  ) \to  T_z(\tilde \Sigma  )$, where $W^{\vec h}(Y)$ is  given by the second fundamental form
associated with $\vec h$:
\beqn\label{eq:12}
\langle W^{\vec h}(Y) , X \rangle = - II_{\vec h}(Y, X) = - \langle \vec h, \nabla_{X}Y \rangle
\eeqn
for all $X \in T_z( \tilde \Sigma)$. Hence, our Jacobi field $J$ satisfies the Coddazzi equation
\beqn\label{eq:13}
[J'(0)]^\top = W^{\vec h} J(0),
\eeqn
where $[\vec \eta]^\top $ denotes the tangential component of $\vec \eta \in T_z(S^n)$.
It follows that
$$
J \in \Gamma_{\sigma_{z, \vec h}, \tilde \Sigma}.
$$

For $J \in \Gamma^0_{\sigma_{z, \vec h}, \tilde \Sigma}$, we further require that
$J'(0) \in T_z(S^n)$. Hence, it follows from \eqref{eq:13} that
\beqn\label{eq:14}
J'(0) = W^{\vec h} J(0)
\eeqn
Because $J(0) \in T_z( \tilde \Sigma)$ and $\dim ( \tilde \Sigma) = k+1$, by \eqref{eq:14} one has
 $\dim [ \Gamma^0_{\sigma_{z, \vec h}, \tilde \Sigma}  ] \le (k + 1)$.

We now prove that $\dim [ \Gamma^0_{\sigma_{z, \vec h}, \tilde \Sigma}  ] \ge (k + 1)$.
Let $\{\vec v_1, ...., \vec v_{k+1}\}$ be an orthogonal basis of $T_z( \tilde \Sigma)$.

It is well-known that, on each curve $s \to y_i(s)$ with $y_i(0) = z$ and $y_i'(0) = \vec v_i$,
there is a unique
a vector field $\{\vec w_i(s) \}$ satisfying
 $\vec w_i(s) \perp T_{y_i(s)}(\tilde \Sigma)$ and
\beqn\label{eq:15}
[\nabla_{y'_i(s)} \vec w(s)]^\perp = 0
\eeqn
with $\vec w_i(0) = \vec h$, where $[\vec \eta]^\perp$ denotes the normal component of $\vec \eta$.
The linear system \eqref{eq:15} has $(n-k-1)$-unknowns and $(n-k-1)$-equations.
Thus, the system \eqref{eq:15} has a unique solution $\vec w_i(s)$ with $\vec w_i(0) = \vec h$.

Let $F_i(t, s) = \E_{y_i(s)}[t \vec w_i(s)]$ and $J_i(t) = \frac{\partial F_i }{\partial s}(t, 0)$. Then
$J_i \in  \Gamma^0_{\sigma_{z, \vec h}, \tilde \Sigma}$ for $i= 1, 2, ..., k+1$. Clearly,
$\{J_1(0), ..., J_{k+1}(0)  \} = \{\vec v_1, ..., \vec v_{k+1} \}
$  are linearly independent. Thus,
$\dim[\Gamma^0_{\sigma_{z, \vec h}, \tilde \Sigma} ]\ge (k+1)$.

(2) If $\tilde \Sigma$ has the first focal point $\sigma(t_0)$ along $\sigma$ with
$0 < t_0 < \frac{\pi}{2}$, then
there must be an orthogonal  Jacobi field $\{ J(t) \}$ along $\sigma$ in $S^n$ with $J(0) \in
T_z( \tilde \Sigma )$ and $J(t_0) = 0$. It is well-known that, in $S^n$, any Jacobi field $\{J(t)\}$
with $J(t_0) = 0$ can be expressed as $J(t) =  \sin (t-t_0) c E(t)$, where $c$ is a non-zero constant and $\{E(t)\}$
is a unit parallel vector field along $\sigma$.

Because $0< t_0 < \frac{\pi}{2}$, we obtain that $J(0) = -(\sin t_0)
c E(0) \ne 0$. Since  $J(0) \in T_z( \tilde \Sigma )$, we see that
$E(0) = - \frac{1}{ \sin t_0 } J(0) \in T_z( \tilde \Sigma )$. It
follows that $J'(0) = \cot (t_0) c E(t) \in \tilde \Sigma  $ and
hence   $J \in \Gamma^0_{\sigma_{z, \vec h}, \tilde \Sigma}$.
\hfill{Q.E.D.}

\bigskip
\noindent {\bf Step 3.}  We will use the Hessian comparison theorem
to show that {\it if $\pi_p: S^{\perp}_p(A, M) \to A'$ is a
Riemannian submersion, then $\tilde \Sigma_{p, q}$ has the first
focal radius $ \ge \frac{\pi}{2}$ in $S^n $, and hence
$\pi_p^{-1}(q)$ is a great circle by Step 1.}

\medskip

We will divide it into two sub-steps:

\smallskip
\noindent
{\it Step 3.1.} Using the Hessian comparison theorem, we  study the decomposition of $T_y(M^n)$ associated with
$\{A, A'\}$ and parallel transports. As an application, we will show that {\it  the covariant derivatives of
horizontal lifting vector fields along each fiber $\tilde \Sigma_{p, q}$ must be vertical}.

\smallskip
\noindent
 {\it Step 3.2.} By Step 3.1, we will construct all Jacobi fields
$J \in \Gamma^0_{\sigma_{z, \vec h}, \tilde \Sigma}$ with {\it
vertical initial derivatives} explicitly. A simple calculation will
show that any non-trivial element $J \in \Gamma^0_{\sigma_{z, \vec
h}, \tilde \Sigma}$ has  the non-vanishing property $J(t) \ne 0 $
for $0 \le t  < \frac{\pi}{2}$. This will complete the proof of
Theorem A.

\bigskip
\noindent {\bf Step 3.1.} { Hessian comparison and the covariant
derivatives of horizontal lifting vector fields along each fiber.}

\bigskip

   We will frequently use the following result.

 \begin{prop}\label{thm:04} (cf. [GG1]) Let $A, A'$ and $M$ be as in Proposition 2.4.
 Suppose that $(p, q) \in A \times A'$. Then

 (1)  Whenever $\dim(A') > 0$, for any unit tangent vector $\vec \eta_0 \in T_q(A')$ and a unit normal vector
$\vec v \in S_q^\perp(A', M)$, the image
$\E_q(\R^2_{\vec \eta_0, \vec v})$ is a totally geodesic immersed 2-sphere of constant sectional curvature
$1$, where $\R^2_{\vec \eta_0 , \vec v} = \text{Span}_{\R}\{\vec \eta_0, \vec v  \}$
is a real 2-dimensional tangent subspace spanned by $\{\vec \eta_0 ,
\vec v  \}$.

(2) Let $\Psi_{A'}: [M - A] \to A'$ be the nearest point projection,
$\Phi_{A'}: ([B^\perp_{ \frac{\pi}{2} }(0_p) - \{ 0_p\}], g_1) \to
A'$ be given by $\Phi_{A'} (z) = \Psi_{A'}(\E_p(z))$ for $z \in
B^\perp_{ \frac{\pi}{2} }(0_p) = \{ z \in T_p(M) | z \perp T_p(A),
|z| < \frac{\pi}{2}\}$ and $z \neq 0$. Suppose that $S^\perp_{ r_0}
(0_p) = \partial B^\perp_{ r_0 } (0_p)  $. Then for $r_0 \in (0,
\frac{\pi}{2}) \to A'$, the map
$$
\Phi_{A'}|_{S^\perp_{ r_0} (0_p)}: S^\perp_{ r_0 }(0_p) \to A'
$$
is a Riemannian submersion up to a constant factor $c = \frac{1}{\sin r_0}$ with respect to
the spherical metric $g_1$ on $ S^\perp_{ r_0} (0_p) \subset B_\pi(0_p)$.

   The similar conclusions hold at $p \in A$ if $\dim(A) > 0$.
\end{prop}

Let us consider the normal bundle of $\tilde \Sigma_{p, q}$ at $z$ with $0 < |z| \le  \frac{\pi}{2} $.

 \begin{definition}  Let $(p, q) \in A \times A'$ be as above. If $z \in \tilde \Sigma_{p, q} \subset B_\pi(0)
\subset T_p(M^n)$ with $0 < |z| \le  \frac{\pi}{2} $ and if $\vec h
\perp T_z( \tilde \Sigma_{p, q}  )$ then the vector $\vec h$ is
called a horizontal vector.

 Similarly, if $\hat z \in M^n$ with $0 < d(p, \hat z) \le  \frac{\pi}{2}$ and
$\hat h \perp T_z(  \Sigma_{p, q}  )$ then the vector $\hat h$ is
called a horizontal vector.

The horizontal subspace at $\hat z$ is denoted by $H_{\hat z}$.
\end{definition}

We will use the Hessian comparison theorem show that the horizontal subspaces is invariant
under the parallel translation along radial geodesics from
$A'$ to $p$. If $c: [a, b]
\to M$ is a curve, we let $\tau_{c(t_1)}^{c(t_2)} $ be the parallel translation along the curve $c$.

\begin{thm}\label{thm:06}  Suppose that $(p, q) \in A \times A'$ and $M^{n}$ are as in Proposition 2.4 and suppose
that  $\sigma: [0, \frac{\pi}{2}] \to M^{n}$ be a geodesic of unit speed from $q$ to $p$. Then
the tangent space $T_{\sigma(t)}(M^n)$ has the following orthogonal decomposition:
$$
T_{\sigma(t)}(M^n) = \tau_{\sigma(0)}^{\sigma(t)}[T_q(A')] \bigoplus
\tau_{\sigma(\frac{\pi}{2})}^{\sigma(t)}[T_q(A)]
\bigoplus T_{\sigma(t)}(\Sigma_{p, q}).
$$
Hence, $\tau_{\sigma(0)}^{\sigma(t)}[T_q(A')]
\bigoplus\tau_{\sigma(\frac{\pi}{2})}^{\sigma(t)}[T_q(A)]  $ is
equal to the horizontal subspace $ \textsl{H}_{\sigma(t)}$ at
$\sigma(t)$ for $t \in [0, \frac{\pi}{2})$.
\end{thm}
{\bf Proof.} By Lemma 3.1 of [GG1], $\sigma'(t) \perp \textsl{H}_{\sigma(t)}$.
We need to show that $ T_{\sigma(t)}( \Sigma_{p, q}  ) \perp
\textsl{H}_{\sigma(t)}$. For this purpose, we use the sharp version of Hessian comparison.

Let $m = \dim A$, $m' = \dim A'$ and $k+1 = \dim(\Sigma_{p, q}  )$. We will also see  that
$\dim M^n = n = m + m' + (k+1)$.

 Let $f(x) = d(x, A')$.  Because $A'$ is totally
geodesic, there are $m'$ Jacobi fields $\{J_1(t), J_2(t), ..., J_{m'}(t) \}$ along $\sigma$ such that
$\{J_1(0),  J_2(0), ..., J_{m'}(0) \}$ is an
orthonormal basis of $T_q(A')$ and $J'_i(0) = 0$ for $i = 1, 2, ..., m'$.

Similarly, if $\dim A > 0$, there are $m$ Jacobi fields
$\{J_{m'+1}(t), J_{m'+2}(t), ..., J_{m'+m}(t) \}$ along $\sigma$
such that $\{J_{m'+1}(0),  J_{m'+2}(0), ..., J_{m'+m}(0) \}$ is an
orthonormal basis of $T_q(A')$ and $J'_{m'+j}(0) = 0$ for $i = 1, 2,
..., m$.

We already knew that the cut-radii of $A'$ and $A$ are equal to the
diameter of $M^{n}$, which is $\frac{\pi}{2}$. Thus, $J_i(t) \neq 0$
for $i = 1, 2, ..., 8$ and $t \in [0, \frac{\pi}{2})$. Recall that
the sectional curvature $\ge 1$, by Berger comparison theorem (or
the 2nd Rauch comparison theorem), we can find a parallel vector
field $\{ E_i(t)\}$ along $\sigma$ such that $J_i(t) = \cos t
E_i(t)$ for $i = 1, 2,...m'$ and $J_{m'+j}(t) = \sin t E_{m'+j}(t)$
for $j = 1, ... m$ if $m = \dim A > 0$. It is clear that
$$
Hess(f)(J, J) = \langle J(t), J'(t) \rangle.
$$
It also is well-known that
the Hessian of distance function $f$ satisfies the so-called Riccati equation:
$$
 \nabla_{ \sigma'(t) }[\text{Hess}(f)] + [ \text{Hess}(f)    ]^2 + R = 0.
$$
More precisely, we let $\{E_i(t) \}_{1 \le i \le n}$ be a parallel
orthonormal base along the geodesic segment $\varphi_v$ with
$E_{n}(t) = \sigma'(t)$, $H_{i, j} (t) = \text{Hess}(f)(E_i(t),
E_j(t))$ and $R_{ij}(t)= \langle R( \sigma(t), E_i(t)) \sigma'(t),
E_j(t) \rangle$, where $R(X, Y) Z = - \nabla_X \nabla_Y Z + \nabla_Y
\nabla_X + \nabla_{[X, Y]}Z$ is the curvature tensor. Thus, we have
$$
H' + H^ 2 + R = 0.
$$

Let $$
 W_{A'}(t) = \{ Y(t) | \quad  H(., Y(t) )|_{ \sigma(t) } = \tan (t) \langle ., Y(t) \rangle \}
 $$
 and
$$
W_{A}(t) = \{ Y(t) | \quad  H(., Y(t) )|_{ \sigma(t) } = \cot (t) \langle ., Y(t) \rangle \}.
$$

 We have shown  that
the eigenspace $\{W_{A'}(t) \}$  is invariant under parallel translation along $\sigma$. Similarly, if
$\dim A > 0$, then $\{W_{A}(t) \}$  is invariant under parallel translation along $\sigma$.

Choose $t_0 = \frac{\pi}{3}$. It is clear $\cot  \frac{\pi}{3} \ne \tan \frac{\pi}{3}$. Thus,
$$
        W_{A'}(t) \perp W_{A}(t)
$$
whenever $\dim A > 0$.

   In what follows, we prove that
   $$
T_{\sigma(t)}(\Sigma_{p, q}) \perp [W_{A'}(t) \bigoplus W_{A}(t)   ].
   $$

We already showed that $E_j( t ) \in W( t)$ for $j = 1, 2, ..., m'$.
Notice that $H_{jj}(t) $ blows up as $t \to 0^+$ for $j > (m'+m)$. If $\{ \lambda_{m+m'+1}(t),
\lambda_{m+m'+2},...,
\lambda_{m+m' + k }(t)\}$ are other eigenvalues of $H$,  then  $\lambda_j(t) \to + \infty$ as $t \to 0$ for
$ j \le (m+m')$. Thus, the corresponding eigenvectors are orthogonal to $W_{A'}(t)$, because eigenvalues are
different.

Similarly, if $\dim A > 0$, we consider $t \to \frac{\pi}{2}$, then $\lambda_j(t) \to + \infty$ as $t \to
\frac{\pi}{2}$ for
$ j \le (m+m')$. For the same reason, the corresponding eigenvectors are orthogonal to $W_{A}(t)$, because
eigenvalues are
different.

Therefore, we proved
$$
H_{ij } (t) = 0
$$
for $i=1, ..., (m+m')$ and $j > (m + m')$.

 Let $\{ (x_1, .., x_{m'})\}$ be a geodesic normal coordinate system of $A'$ at $q$ given by
 $G: \R^{m'} \to A'$ with $G(x_1, ..., x_{m'}) = \E_q(\sum_1^{m'} x_i
 E_i(0)   )$. Recall that $\dim \{ [  T_q(A') ]^\perp \} = m + k + 1$. Thus there exists an
  orthonomal basis $\{E_{m'+1}, ..., E_{n-1}, E_{n}\}$  of
 $ [  T_q(A') ]^\perp$
 such
 that $E_{n} = \sigma'(0)$. Let $\vec \theta  = (\theta_{m'+1}, ..., \theta_{n})$
 with $|\vec \theta| \le 1$. Then
$\Psi: B_1(0) \to S^{n-m'-1} = S^\perp_q(A', M^{n})$ given by
$$\Psi(\theta_{m'+1}, ..., \theta_{n-1}   ) =
\sum_{j = m'+1}^{n} \theta_j E_j + \sqrt{1- |\vec \theta|^2} \, \, \sigma'(0)
$$
gives rise to a local  coordinate system  of  $S^{n - m' -1} = S^\perp_q(A', M^{n})$ around $\sigma'(0)$.
 Using the parallel transport $\tau_{G(0)}^{G(x)}$ from $q = G(0)$ to $G(x)$ we have a local
 coordinate system
 given by $ (\theta_{m'+1}, ..., \theta_{n-1}) \to   \tau_{G(0)}^{G(x)}(\Psi (\vec \theta )  ) $
 for $S^\perp_{G(x)}(A', M^{n})$.
 Therefore, $\{ (x_1, ..., x_{m'}; \theta_{m'+1}, ..., \theta_{n-1}, t)\}$ gives rise to a
 local coordinate for normal bundle of $A'$ in $M^{n}$ near $(x, t\sigma'(0))$. In fact, the
$F( x_1, ... x_{m'};  \theta_{m'+1}, ...,\theta_{n-1}, t) = \text{Exp}_{G(x)}
[ t \tau_{G(0)}^{G(x)}(\Psi (\vec \theta )  )   ]$ does the job.
Finally we let
$C_{ji}(t) =
\langle \frac{\partial F}{\partial \theta_j  }, E_i \rangle|_{\sigma (t)}$.
It is well-known that $H(t) = C'(t) [ C(t) ]^{-1}$.
We already showed  that $ W_{A'}(0) = T_q( A') $ and $\{ W(t) \}$ is
parallel along $\sigma$.
Using $C_{ji}( 0 ) = 0$   and the
fact $H_{ij } (t) = 0$ for $i=1, ..., m'$ and $j > (m' + 1)$, by the integration of $C'(t) = H(t) C(t)$
from $\frac{\pi}{2}$ to $t$ we conclude that
$$
C_{ij } (t) = 0
$$
for $i=1, ..., m'$ and $j > (m'+1)$. Thus,
we see that $\frac{\partial F}{\partial \theta_j  } \in [ W_{A'}(t) ]^\bot$
for $j > (m' + 1)$.

   Therefore, both tangential subspace $T_{\sigma(t)}(\Sigma_{p, q}     ) \bigoplus
   W_A(t) $ and  sub-space $W_{A'}(t)  $ at $\sigma(t)$ are
   invariant under parallel translation along $\sigma$.
It follows that
   $$
   T_{\sigma(t)}(\Sigma_{p, q}     ) \perp W_{A'}(t) .
   $$

   For the same reason, if $\dim A > 0$,  one has
    $$T_{\sigma(t)}(\Sigma_{p, q}     ) \perp W_{A}(t)
    $$
    as well. We already proved $W_{A}(t) \bot W_{A'}(t)$. This completes the proof.
\hfill{Q.E.D.}

Theorem \ref{thm:06} indicates that there is a  non-trivial relation
between the exponential map $\E_A$ along the normal bundle of $A$
and the exponential map $\E_{A'}$ along the normal bundle of $A'$.
As an application of Theorem \ref{thm:06}, we draw some conclusions.

 \begin{cor} Let $(p, q) \in A \times A'$, $A$ and $A'$ be as in Proposition \ref{thm:02}
 and $\dim A' > 0$. Suppose that $\vec \eta \in T_q(A')$, $\hat z \in \Sigma_{p, q}$ with
 $0< d(\hat z, A') < \frac{\pi}{2} $,  $\hat h_{\eta}(\hat z)$ is the parallel transport of $\vec \eta$ along
 the unique length-minimizing geodesic segment from $q$ to $\hat z$,
 $z = (\E_p)^{-1}(\hat z) \in \tilde \Sigma_{p, q} $ and
 $$
\vec h_{\eta}(z) = [(\E_p)^{-1}_*]|_z \hat h_{\eta}(\hat z).
 $$
Then the horizontal lifting vector field $\{ \vec h_{\eta}\}_{z \in \tilde \Sigma_{p, q}}$ of $\eta $ has the property
\beqn\label{eq:21}
\nabla_X \vec h_{\eta} \in T_z( \tilde \Sigma_{p, q} )
\eeqn
for all $X \in T_z( \tilde \Sigma_{p, q} )$.
\end{cor}
{\bf Proof.}  We first consider the case $X = \nabla r $, where $
r(z) = |z| = d(0_p, z)$. By our assumption, there is a unique
geodesic segment of unit speed  from $q$ to $\hat z$, say
$\sigma_{q, \hat z}$. Let $\vec v = \sigma_{q, \hat z}'(0)$. By
Proposition \ref{thm:04}, if we let $\R^2_{\vec \eta, \vec v} $ be
the subspace spanned by $\{ \vec \eta, \vec v\}$, then $\hat S^2 =
\E_q( \R^2_{\vec \eta, \vec v} )$ is a totally geodesic immersed
2-sphere $S^2_{\vec \eta, \vec v}$ of constant curvature 1, which
passes both $p$ and $q$. It follows that,  on the unit 2-sphere
$S^2_{\vec \eta, \vec v}$, on has \beqn\label{eq:22} \nabla_{\nabla
r } \vec h_{\eta}|_z = 0 \eeqn

  We now consider the remaining case $X \in  T_z( \tilde \Sigma_{p, q} )$ but $X \bot \nabla r$. Let
 $$
 B_{0_p, \pi}^\perp =  \{ z \in T_p(M) \,  |  \, z \perp T_p(A), |z| \le \pi \}.
$$
In terms of the spherical metric $g_1$, the sub-manifold $(B_{p,
\pi}^\perp, g_1) $ is a totally geodesic $(n-m)$-dimensional sphere
$S^{n-m}$, where $m= \dim A$, $m'= \dim A'$, $k+1 = \dim \Sigma_{p,
q}$ and $n = \dim M = m+ m' + k + 1$.

Let $\Psi_{A'}: [M-A] \to A'$ be the nearest point projection,
$S_{p, r}^\perp = \{ z \in T_p(M) \,  |  \, z \perp T_p(A), |z| = r
\}$. By Proposition \ref{thm:04} and Theorem \ref{thm:06}, in terms
of the spherical metric $g_1$ on $B_\pi(0_p)$, the map
$$
\begin{array}{rrcl}
\tilde \Psi_{A'}: &  S_{p, r}^\perp & \rightarrow & A' \\
                          & z    & \rightarrow & \Psi_{A'}[\mathrm{Exp}_p(z)]
\end{array}
$$
is a Riemannian submersion up to a constant factor $\frac{1}{\sin r}  $.
Consequently, if $\{\vec \eta_1, ...., \vec \eta_{m'}\}$ is an orthonormal basis of
$T_q(A')$, then by Theorem \ref{thm:06} and its proof, the set of vectors
$$
\{\vec h_{\vec \eta_1}, ...., \vec h_{ \vec \eta_{m'}}\}
$$
form a basis of the normal bundle  $N(\tilde \Sigma_{p, q},  B_{0_p, \pi}^\perp )$ of $ \tilde \Sigma_{p, q}$
at $z$ in $S^{n-m} = ( B_{0_p, \pi}^\perp, g_1 )$.

Let $\{x_1, ...., x_{m'}\}$ be the geodesic normal coordinate of
$A'$ at $q$. We choose $\vec \eta_i = \frac{\partial  }{\partial
x_i}$ at $0_q$; i.e., we use the map $(x_1, ..., x_{m'}) \to
\E_q(x_1 \vec \eta_1 + ... x_{m'} \vec \eta_{m'}   )$ as the
geodesic  coordinate system of $A'$ at $q$.

  Suppose that $G(x) = \E_q(x_1 \vec \eta_1 + ... x_{m'} \vec \eta_{m'}   )$
  and recall that $\vec v = \sigma'_{q, \hat z}(0)$.
  Let us consider the Fermi coordinate system
   (the exponential map) along $A'$:
$$ \digamma (x, \rho  \vec v) = \E_{G(x) }[\tau_{q}^{G(x)} \rho \vec v  ].
   $$
By the proof of Theorem \ref{thm:06}, we have
\beqn\label{eq:23}
\vec h_{\vec \eta_i} = \frac{1}{\sin |z|} [\frac{\partial (\E_p^{-1} \circ \digamma) }
{\partial x_i}]|_{(0_q, \rho_0\vec v)},
\eeqn
where $$\rho_0 = \frac{\pi}{2} - |z|.$$

For simplicity, we denote $\E_p^{-1} \circ \digamma  $ by $  \tilde F$.
We now choose $X = \frac{\partial \tilde F }{\partial v_i}|_{(0_q, \rho_0\vec v)}$ for
$\vec v = (v_1, ..., v_k) \in [ T_q(A') ]^\perp$ and
$ |\vec v | = 1$, where we only allow $\vec v \in S^\perp_q(A', M)$.
It is easy to see (cf. [CE, page2]) that, if $[X, Y] = 0 = [Y, Z] = [X, Z]$, then
$$
\langle  \nabla_XY , Z \rangle = X\langle  Y , Z \rangle +  Y\langle  X , Z \rangle -
Z\langle  X , Y \rangle.
$$
Therefore, setting $ \rho_0 = \frac{\pi}{2} - |z| $ and by a direct calculation, one has that, if
$X = \frac{\partial \tilde F }{\partial v_i}|_{(0_q, \rho_0\vec v)}$ then
$$
\langle \nabla_X \vec h_{\vec \eta_i}, \vec h_{\vec \eta_j}\rangle = 0 + 0 - 0 = 0
$$
for all $i, j = 1, ..., m'$. This completes the proof.
\hfill{Q.E.D.}
\bigskip

A direct consequence of the above corollary is the following result.

\begin{cor}\label{thm:09} Let $(p, q) \in A \times A'$, $A$ and $A'$ be as in Proposition \ref{thm:02}
and $\dim A' > 0$. $\vec \eta \in T_q(A')$, $\hat z \in \Sigma_{p,
q}$ with $0< d(\hat z, A') < \frac{\pi}{2} $,  $\hat h_{\eta}(\hat
z)$ is the parallel transport of $\vec \eta$ along the unique
length-minimizing geodesic segment from $q$ to $\hat z$, $z =
(\E_p)^{-1}(\hat z) \in \tilde \Sigma_{p, q} $,
 $$
\vec h_{\eta}(z) = [(\E_p)^{-1}_*]|_z \hat h_{\eta}(\hat z)
$$
 and
 \beqn\label{eq:24}
F_{\vec \eta}(t, z) = \E^{S^n}_z[t \vec h_{\eta}(z)  ].
 \eeqn
Then the corresponding Jacobi fields
\beqn\label{eq:25}
 J_i(t) =   \frac{\partial F_{\vec \eta} }{\partial z_i}(t, z)
\eeqn
along the geodesic $\{F_{\vec \eta}(., z)\}$
has the property
\beqn\label{eq:26}
J_i'(0) \in T_z(\tilde  \Sigma_{p, q}  )
\eeqn
for $i=1, ..., k+1$, where $\{z_1, ...., z_{k+1}) \}$ is any local coordinate system of
$\tilde \Sigma_{p, q}$ around $z$.
\end{cor}

\bigskip
\bigskip
\noindent {\bf Step 3.2.}  Proof of Theorem A.

\bigskip

We recall some elementary facts about the geodesic triangles in a unit 2-sphere $S^2$, which are isometrically
immersed in $M$, see Proposition 3.3 (1) above.

\begin{lemma}\label{thm:10} Let $ \hat z \in \Sigma_{p, q}$ with $0 < r_0 = d(\hat z, A) < \frac{\pi}{2}$
and $S^2 =S^2_{ q, (z, \hat h)}$ be a totally
geodesic immersed 2-sphere in $M$  given  by $S^2_{ q, (z, \hat h)}=
\E_{\hat z}( \R^2_{q, (z, \hat h)} )$, where $\hat h$ is a unit
horizontal vector and $\R^2_{q, (z, \hat h)} = Span\{\hat h,
\E_{\hat z}^{-1}(q)\}$ described as in Proposition 3.3(1).

(1) If $\varphi_{\hat h} (t) = \E_{\hat z}( t \hat h)$ for some unit horizontal
vector at $\hat z$ with
$0 < r_0 = d(\hat z, A) < \frac{\pi}{2}$, then
$\varphi_{\hat h}( \frac{\pi}{2} ) \in A'$;

(2) For $0 < t < \frac{\pi}{2}$, the distance function satisfies
 $r(t) = d( \varphi_{\hat h} (t) , A) = \arccos [ \cos r_0 \cos t] $;
 Consequently, $d(\varphi_{\hat h} (t), A')
 = \arcsin [ \cos r_0 \cos t   ]$, where $r_0 = d( \hat z, A)$.

 (3) Let $\Psi_{A'}: [M-A] \to A'$ be the nearest point projection and $\ell(t) $ be the length
 of $\Psi_{A'}[\varphi_{\hat h}([0, t])  ]$. Then
 $$
\ell (t)  = \arccos [\frac{ \sin r_0 \cos t   }{ \sqrt{1 - (\cos r_0 \cos t)^2  } }     ].
 $$

 (4) The vector $[\varphi_{\hat h} '(t) - \langle \varphi_{\hat h} '(t), \nabla r \rangle \nabla r ]$
 remains to be horizontal, where $r(x) = d(A, x)$.

\end{lemma}

The lemma above can be proved by the law of cosine in $S^2$, see [Pe, page 314].

Finally, we can now show that $\tilde \Sigma_{p, q}$ has focal radius $\ge \frac{\pi}{2}$ in $S^n$.

\begin{lemma}\label{thm:11} Let $z \in \tilde \Sigma_{p, q} \subset S^n$ and $J_i(t)$ be
as in Corollary \ref{thm:09} above. Then
$$
J_i(t) \neq 0
$$
for all $t \in (0, \frac{\pi}{2}  )$. Consequently, $\tilde \Sigma_{p, q}$ has
focal radius $\ge \frac{\pi}{2}$ in $S^n$.
\end{lemma}
{\bf Proof.} We choose a special local coordinate system of $ \tilde \Sigma_{p, q} $ at $z$ as follows.
By Corollary \ref{thm:09}, $J_i'(0) \in T_z( \tilde \Sigma_{p, q} )$ for all $i = 1, ..., k+1$. We can choose
$(k+1)$-principal directions $\{ e_1, ..., e_{k + 1}\}  $
of the Weingart map $W^{\vec h}: X \to (\nabla_X \vec h(z)   )^\top = \nabla_X \vec h(z)   $ for all
$X \in T_z(M)$, where
$$
\langle W^{\vec h} X, Y \rangle = \langle \nabla_X \vec h(z), Y \rangle
$$
for all $X, Y  \in T_z(M)$ and $(\vec w)^\top$ is the tangential component of $\vec w  $.

It was proved $ \nabla r|_z = \frac{z}{|z|}$ and $\vec h(z)$ span a
totally geodesic 2-sphere of constant curvature 1, see Proposition
3.3 (1) above. Thus, $\nabla r|_z$ is an eigenvector of $W^{\vec
h}$.  We choose $e_{k +1 } = \nabla r|_z$. Furthermore, in $S^2$,
the corresponding Jacobi field can be written as $J_{k+1}(t) = (\cos
t) E_{k +1} (t)$, where $\{ E(t) \}$ is a parallel vector along
$\sigma_{z, \vec h}(t) = \E_z( t \vec h)$ with $E(0) =  \nabla
r|_z$. Hence, $J_{k+1}(t) \neq 0$ for all $t \in (0,
\frac{\pi}{2})$.

\bigskip

We now consider the remaining $\{J_1, ..., J_k \}$. Let
$$
   \upsilon_i(s) = |z|[ (\cos s) \frac{z}{|z|} + (\sin s) e_i]
$$
and
$$
F_i(t, s) = E^{S^n}_{\upsilon_i(s)}(t \vec h_{\vec \eta}( \upsilon_i(s)  ))
$$
for $i = 1,..., k$. Finally, we set
$$
J_i(t) = \frac{\partial F }{\partial s}(t, 0)
$$
for $i = 1, ..., k$.

In order to prove that $J_i(t) \neq 0$ for $t \in (0, \frac{\pi}{2})$, we use
$$
     \hat F_i(t, s) = \E_p^M[ (\E_p^{S^n}) ^{-1} (F_i(t, s)    )]
$$
for $i = 1, ..., k$.
By Lemma \ref{thm:10}, one has
\beqn\label{eq:27}
0 < \rho(t) = d(A', \hat F_i(t, s)) = \arcsin [ (\cos |z|) \cos t] < \frac{\pi}{2}
\eeqn
for $ 0 \le t < \frac{\pi}{2}$.
Let $G = \E_p^M[ (\E_p^{S^n}) ^{-1}]$ and
$$
\hat J_i(t) = \frac{\partial \hat F_i }{\partial s}(t, 0) = G_* J_i(t).
$$
Because $G$ is a local diffeomorphism at all $x \in B_{\frac{\pi}{2}}(0_p)$ with $0 < |x| < \frac{\pi}{2}$,
using \eqref{eq:27} one concludes the following is true: {\it ``$J_i(t) \neq 0$  holds for
$ t \in (0,  \frac{\pi}{2})$ if and only
if $\hat J_i(t) \neq 0$ holds for  $ t \in (0,  \frac{\pi}{2})$"}.

\medskip

It remains to verify that $\hat J_i(t) \neq 0$ for $ t \in (0,  \frac{\pi}{2})$. For this purpose, we
express $\hat J_i(t)$ in terms of the Fermi coordinates along $A'$ instead. In terms of the
Fermi coordinates along $A'$, we will clearly see that $\hat J_i(t) \neq 0$ for $ t \in (0,  \frac{\pi}{2})$.
The detail for the new expressions of $\hat J_i(t)$ and $\hat F_i(t, s)$ can be given as follows:

Notice that $\{\hat h_{\vec \eta}(\upsilon_i(s) ), \E_{\upsilon_i(s)}^{-1}(q)\}$ span a totally geodesic
immersed 2-sphere $S^2_{\upsilon_i(s), \hat h_{\vec \eta}}$ of constant curvature $1$. Such a 2-sphere
$S^2_{\upsilon_i(s), \hat h_{\vec \eta}}$ passes through the geodesic $\hat \sigma_{\vec \eta}(\ell) =\E_q
(\ell \vec \eta)$.
Let
\beqn\label{eq:28}
\vec \psi_i(s) = \E^{-1}_q[ \hat F_i(0, s)  ],
\eeqn
$\tau_{q}^{\hat \sigma(\ell)}$ be the parallel translation along $\hat \sigma_{\vec \eta}$ and
let $\Psi_{A'}: [M- A] \to A'$ be the nearest point projection. Then,  by Lemma \ref{thm:10}, one has
 $$
\ell (t) = d(\Psi_{A'}( \varphi_{\hat h} (t), q  ) = \arccos
[\frac{ \sin r_0 \cos t   }{ \sqrt{1 - (\cos r_0 \cos t)^2  } }
].
 $$
 A direct calculation shows that if $q(t) = \hat \sigma_{\vec \eta}(\ell (t))$ then
 \beqn\label{eq:29}
\hat F_i(t, s) = \E_{q(t)}[ \tau_{q}^{q(t)}
\rho(t) \vec \psi_i(s)  ]
\eeqn
for $i = 1, ..., k$. It follows that
\beqn\label{eq:30}
\hat J_i(t) = [\E_{q(t)}]_* [  \tau_{q}^{q(t)}( \rho(t) \vec \psi_i'(0)) ]
\eeqn

We already proved that $0< \rho(t) < \frac{\pi}{2}$. Recall that $A'$ is totally geodesic and
\beqn\label{eq:30}
[\tau_{q}^{q(t)}( \rho(t) \vec \psi_i'(0))] \perp T_{q(t)}(A')
\eeqn
for all $t$. Recall that the parallel transport $\tau_{q}^{q(t)}: T_q(M) \to T_{q(t)}(M)$ is an isometry.
Since the cut radius of $A'$ is equal to $\frac{\pi}{2}$, it follows equations \eqref{eq:27} -\eqref{eq:30}
that
$$
\hat J_i(t) = [\E_{q(t)}]_* [  \tau_{q}^{q(t)} \rho(t) \vec \psi_i'(0) ] \neq 0
$$
for $i = 1, ..., k$ and $t \in (0, \frac{\pi}{2} )$,  as long as
$\vec \psi_i'(0) \neq 0$ and $\rho(t) \neq 0$. Recall that $J_i(0)
\neq 0$ and $\rho(t) \neq 0$ for $t \in (0, \frac{\pi}{2} )$. Hence
$\vec \psi_i'(0) \neq 0$ and $\hat J_i(t) \neq 0$ holds  for $i = 1,
..., k$ and $t \in (0, \frac{\pi}{2} )$. This completes the proof.
\hfill{Q.E.D.}

\bigskip
\noindent{\it The end of the proof of Theorem A.}\
   By Steps 1-3 above, we proved that $\pi^{-1}_p(q)$ is a
great circle for each $(p, q) \in A \times A'$.
   Furthermore, it follows from Proposition 3.1(2) that $\mathrm{Diam}(A') = \frac{\pi}{2}$.
   We can also choose a point $y \in A'$ with $d(y, q) =\frac{\pi}{2}$. Using Proposition 3.1(2) again,
   we see that $\Inj_M(y)=\frac{\pi}{2}$. By replacing $p$ by $y$ if needed, we may always assume that
$\Inj_M(p)= \frac{\pi}{2}$ and $\dim A = 0$.
Hence, by [Ran], $\pi_p: S_p(M) \to A'$ is isometric to the classical Hopf fibration and $M$ is
isometric to one of $\{ \C P^{\frac n2}, \HH P^{\frac n4}, \C aP^2 \}$.

\medskip

    Professor Grove kindly pointed out that ``if $\pi_p: S_p(M) \to A'$ is a great circle fibration then
    one can show that $M^n$ is isometric to one of $\{ \C P^{\frac n2}, \HH P^{\frac n4}, \C aP^2 \}$ directly
    without using [Ran]." The following argument is an outline of a direct proof inspired by Professor
    Grove, but authors are
    responsible for all possible errors.

Let $y \in A'$ with $d(y, q) =\frac{\pi}{2}$ be as above and $M'$ be
the convex hull of $\{y\} \cup \Sigma_{p, q}$ in $M$. Then, by the
$\pi$-convexity of $S_{\frac {\pi}{2}}(y)$ described in [GG1], one
has $\pi_y: S_y(M') \to \Sigma_{p, q}$ a Riemannian submersion as
well. Steps 1-3 above implies that $\pi_y: S_y(M') \to \Sigma_{p,
q}$ is a great circle fibration. For any great circle fibration
$\pi_y: S_y(M') \to \Sigma_{p, q}$, using O'Neill formula, one can
easily show that $\Sigma_{p, q}$ is isometric to a round sphere of
constant curvature $4$. Thus, each fiber $\Sigma_{p, q}$ is
isometric to $S^{k+1}$ up to a factor $\frac 12$.

     The metric of $(M^n, g)$ can now be explicitly expressed as follows.

   Recall that $M = \cup_{q \in A'} \Sigma_{p, q}$. For each $\hat z \in
    M^n$ and $\xi \in T_{\hat z}(M)$ with $r(\hat z) = d(p, \hat z)$, we let $\xi^H $ denote the horizontal
    component of $\xi$ and we let $\xi^v $ denote the vertical
    component of $\xi$.

    Since each $\Sigma_{p, q}$ is isometric to $S^{k+1}$ up to
    a factor $\frac 12$ and $\pi_p: S_p(M) \to A'$ is a great circle fibration, we have
\beqn\label{eq:31}
|\xi|_g^2 = (\sin r)^2 | \xi^H |^2  + [\frac 12 \sin (2r)   ]^2  | \xi^v |^2.
\eeqn
Using \eqref{eq:31} and an induction method on $\frac{\dim M}{k}$, one can show  that $(M, g)$ is
isometric to one of $\{ \C P^{\frac n2}, \HH P^{\frac n4}, \C aP^2 \}$.
\hfill{Q.E.D.}

\medskip

\noindent{\bf Acknowledgement.} In the earlier version of the manuscript, authors mistakenly viewed
the extrinsic cut-radius of a sub-manifold as its intrinsic injectivity radius. We are very grateful
to Professor Karsten Grove for pointing out the mistake. Professor Grove also kindly suggested the new title for our
paper. Dr. Vitali Kapovitch also provided some suggestions on the writing.

\vspace{1 cm} \centerline{\bf \large References}
{\small\renewcommand{\baselinestretch}{2}

\noindent [Ca] \ E. Calabi, {Hopf's maximum principle
with an application to Riemannian geometry }, {\it
Duke Math. J.}, Vol. 18 (1957), 45-56.

\noindent [CE] \  J. Cheeger and D. Ebin, {Comparison Theorems in Riemannian Geometry},
 North-Holland Publishing Company,  New York,
1975.

\noindent [GG1] \ D. Gromoll and K. Grove, {A generalization of
Berger's rigidity theorem for positively curved manifolds}, {\it
Ann. Scient. Ec. Norm. Sup.}, Vol. 20 (1987), 227-239.

\noindent [GG2] \ D. Gromoll and K. Grove, {The low-dimensional
metric foliations of Euclidean spheres}, {\it J. Diff. Geom.}, Vol.
28 (1988), 143-156.

\noindent [Grv] \ K. Grove, {Privative communications. }

\noindent [GS] \  K. Grove and K. Shiohama, {A generalized sphere theorem}, {\it Ann. of Math.}, Vol. 106 (1977),
371-376.

\noindent [Pe] \  P. Petersen, {Riemannian Geometry}, {\it Graduate Texts in Mathematics}, vol 171, Springer, New York,
1997.

\noindent [Ran] \ A. Ranjan, {Riemannian submersion of spheres with
totally geodesic fibres}, {\it Osaka J. Math.}, Vol. 22 (1985),
243-260.

\noindent [Wil] \ B. Wilking, {Index parity of closed geodesics and
rigidity of Hopf fibration}, {\it Invent. Math.}, Vol. 144 (2001),
281-295.

\end{document}